\documentclass{article}
\usepackage{amsmath,amssymb}
\usepackage{amsmath,amssymb}

\usepackage{mathrsfs} 

\def\real{I\kern -0,37 em R}
\def\nat{I\kern -0,37 em N}
\def\comp{l\kern -0,50 em C}
\def\que{l\kern -0,50 em Q}
\def\ze{Z \kern -0,50 em Z}
\newcommand{\R}{{I\kern -0,37 em R}}
\newcommand{\K}{{I\kern -0,37 em K}}
\newcommand{\OR}{{\overline{{I\kern -0,37 em R}}}}
\newcommand{\OK}{{\overline{{I\kern -0,37 em K}}}}
\newcommand{\I}{{\bf{I}}}

\newcommand{\qed}{\enspace\vrule  height6pt  width4pt  depth2pt}

\newenvironment{proof}{\par\noindent{\bf Proof.}}{$\qed$\par\bigskip}

\title{Natural Topologies on Colombeau Algebras\footnote{ 
2000 Mathematics Subject Classification: Primary 46F30 Secondary 46T20.\protect\\ Keywords and phrases: Colombeau algebra, generalized function, sharp topologies.}}

\author
{J. Aragona   \and R. Fernandez  \and S. O. Juriaans\footnote{Research partially sponsored by CNPq-Brazil Proc. 300652/95-0}}

\date{} 

\newtheorem{theorem}{Theorem}[section]
\newtheorem{definition}[theorem]{Definition}
\newtheorem{lemma}[theorem]{Lemma}
\newtheorem{corollary}[theorem]{Corollary}
\newtheorem{proposition}[theorem]{Proposition}
\newtheorem{remark}[theorem]{Remark}
\newtheorem{example}[theorem]{Example}

\begin{document}
\maketitle

\begin{abstract}
We define intrinsic, natural and metrizable topologies ${\mathscr T}_{\Omega}$, ${\mathscr T}$, ${\mathscr T}_{s,\Omega}$ and ${\mathscr T}_s$ in ${\mathscr G}(\Omega)$, $\OK$, ${\mathscr G}_s(\Omega)$ and $\OK_s$ respectively. The topology ${\mathscr T}_{\Omega}$ induces ${\mathscr T}$, ${\mathscr T}_{s,\Omega}$ and ${\mathscr T}_s$. The topologies ${\mathscr T}_{s,\Omega}$ and ${\mathscr T}_s$ coincide with the Scarpalezos sharp topologies.
\end{abstract}

\section*{Introduction}

It is well known that the Colombeau's theory was developped aiming to solve non-linear problems of PDEs and ODEs. In a more specific way, Colombeau's full algebra (denoted here  as ${\mathscr G}({\Omega})$, where ${\Omega}$ is an open subset of $\R^m$)  was introduced to be the universe that contains a large amount of  solutions of PDEs and ODEs defined on ${\Omega}$, which are not solvable   classicaly [of course, ${\mathscr G}({\Omega})$ has many other inhabitants, but they do not matter here]. The mere fact that ${\mathscr G}({\Omega})$ is an algebra that contains canonically a copy of ${\mathscr D}'({\Omega})$ has represented a great advance because it made possible the multiplication of distributions without any kind of restriction. But this progress, despite its importance, was merely algebraic, letting aside all the resources of Functional Analysis. Since the machinery of classical Functional Analysis has shown to be prolific, it was expected that, soon or late, it would be possible to define a suitable topology [i.e., compatible with the algebraic structure], in order to have, in this context, a complete set of algebraic and topological tools. The first and most important step in this direction was done in 1995 by D. Scarpalezos in [S1] and [S2], who defined the ``sharp topology'', on the simplified version of Colombeau's algebra. In this work we introduce the topology $\mathscr T_{\Omega}$ on the full algebra ${\mathscr G}({\Omega})$ and prove that the topologies induced by $\mathscr T_{\Omega}$ on ${\mathscr G}_s(\Omega)$ and $\OK_s$ [ring of Colombeau's generalized scalars (simplified version)] coincide with Scarpalezos's sharp topologies. The topology $\mathscr T$ induced by $\mathscr T_{\Omega}$ on $\OK$ [ring of Colombeau's generalized scalars (full version)] made it a topological ring so that ${\mathscr G}(\Omega)$ is a topological algebra on $\OK$. Our starting point differs from the one of Scarpalezos. In fact, using a partial order relation $\leq$ on $\OR$ (defined in [A-J-O-S] and [G-K-O-S], see Def.2.2) we have been able to define a notion of ``generalized semi-norms'' on ${\mathscr G}(\Omega)$ from which the definition of $\mathscr T_{\Omega}$ follows naturally. This way of working suggests an interesting parallel between Colombeau's algebras and the theory of locally convex spaces. 

It is pertinent to say that an important part of the motivation for this work was the resolution of a boundary value problem for a nonlinear PDE that will appear elsewhere. The resolution of this problem involves the completeness of $\mathscr T_{\Omega}$, the definition of a topology $\mathscr T_{\overline{\Omega},b}$ on ${\mathscr G}(\overline{\Omega})$ (where $\Omega$ is a bounded open subset of $\R^m$), the completeness of  $\mathscr T_{\overline{\Omega},b}$, besides many other questions. Of course, the addition of all this material  would make this paper excessively large, which makes impossible to include these results here.

Next, we present the general ideas on which this work is based. Let $\Omega$ be a non-void open subset of $\real^m$ and  ${\mathscr G}(\Omega)$ as in [A-B, Not. 2.1.1]. If  $(\Omega_l)_{l\in\nat}$ is an exhaustive sequence of open subsets of $\Omega$ then it is well known that the natural locally convex topology $t_{\Omega}$ on ${\mathscr C}^{\infty}(\Omega)$ can be defined by the fundamental system of semi-norms $|||\,.\,|||_{\beta , l}$ given by $|||u|||_{\beta , l}:=\sup_{x\in\Omega_l}|\partial^{\beta}u(x)|$ for   $u\in{\mathscr C}^{\infty}(\Omega)\,,\,\beta\in\nat^m$ and $l\in\nat$. \typeout{The definition of our full sharp topology ${\mathscr T}_{\Omega}$ on ${\mathscr G}(\Omega)$ is the transport to ${\mathscr G}(\Omega)$ of the above definition of $t_{\Omega}$.} 

To define our topology ${\mathscr T}_{\Omega}$ we proceed as follows. First, for every $\beta\in\nat^m$ and $l\in\nat$, we define the ``generalized semi-norm'' $||\,.||_{\beta , l}$ on ${\mathscr G}(\Omega)$ by

\vskip-0.4cm
$$||f||_{\beta , l}:=cl(\varphi\in A_0\longmapsto |||\widehat f(\varphi,\,.\,)|||_{\beta , l}\in\real)\, , \,\,\forall\,\,f\in{\mathscr G}(\Omega)\,,$$

\vskip-0.1cm
\noindent
where $\widehat f$ is any representative of $f$. Next we  define \typeout{the $||\,.||_{\beta , l}$-balls (of center $0$) getting} a fundamental system of $0$-neighborhoods of ${\mathscr G}(\Omega)$, determining ${\mathscr T}_{\Omega}$. To this end  we need to introduce a partial order relation $\leq$ on $\OR$ which \typeout{allows us   to write things as $||f||_{\beta , l}\leq a $ ($a\in\OR$). We then introduce a partial order relation $\leq$ on $\OR$ which} generalizes the partial order $\leq$ on $\OR_s$ introduced in   [A-J-O-S] and [G-K-O-S]. \typeout{The fact that the partial order $\leq$ is not a linear one, introduces a new difficulty to lead with $||\,.||_{\beta , l}$-balls. The solution of this problem was to use as} To get  a compatible topology we allow as radius of our balls \typeout{ ``radius'' of the $||\,.||_{\beta , l}$-balls, the} elements of the set ${\cal Q}^{\bullet}:=\{\alpha^{\bullet}_r\,|\,r\in\real\}$, which is  a linearly ordered subset of $\OR$ (see Example \ref{2.3}) formed by $q$-positive units of $\OR$ (i.e. ${\cal Q}^{\bullet}\subset Inv(\OR)\cap\OR_+$). The consideration of the generalized scalars $\alpha^{\bullet}_r$ is the natural generalization of an idea of [A-J]. \typeout{, that is, the use of the elements $\alpha_r\in\OR_s$ represented by $\varepsilon\longmapsto\varepsilon^r$, whose image by the natural immersion $\OR_s \hookrightarrow\OR$ is precisely $\alpha^{\bullet}_r$.}

The same    ideas, together with the obvious definition of the absolute value for elements of $\OK$, lead  to ${\mathscr T}_{\Omega}$ \typeout{to allow write inequalities of the kind $|x|\leq\alpha^{\bullet}_r$ ($r\in\real$) and with it we can} and allow us to introduce   a natural topology ${\mathscr T}$ on $\OK$. It is then proved   that ${\mathscr T}_{\Omega}$ and ${\mathscr T}$ induce on the simplified algebras ${\mathscr G}_s(\Omega)$ and $\OK_s$ the well-known sharp topologies (see \cite{S1}, \cite{S2} and [A-J]).

The best way to endow  an algebra (or ring) $E$ with a topology compatible with its algebraic structure  is to give a filter basis ${\mathscr B}$ on $E$ verifying a   set of axioms which essentially guarantees the continuity of the algebraic operations of $E$. This set of axioms is slightly different than the one   usually appearing in text  books on TVS. This is because a $\OK$-algebra is not  a $\K$-TVS. So, to make the paper self contained, we present, in  the section 1, some   results  about these axioms.


\section{Some basic facts on topological algebras}

 In this  section we present three well know results which guarantee that the topology determined by a given     filter basis $\,{\mathscr B}\,$ on a ring (or module or algebra) is compatible with the ring (or module or algebra) structure. The proofs of these three results follow easily from [B, Ch.3, $\S$ 1, n$^{\rm o}$ 2(pg 12); $\S$ 6, n$^{\rm o}$ 3(pg 75-76)] and are omitted. 
 
 In what follows the word ``ring'' means ``commutative ring with unit'' and the word ``algebra'' means ``commutative algebra with unit''.
 
 \begin{proposition}\label {1.1} Let $A$ be a ring and ${\mathscr B}$ a filter basis on $A$ verifying the following conditions:
 
\begin{description}
	\item $(GA'_I)$ for every $U\in {\mathscr B}$ there exists $V\in {\mathscr B}$ such that $V+V\subset U$;
 
 \item $(GA'_{II})$ for every $U\in {\mathscr B}$ there exists $V\in {\mathscr B}$ such that $-V\subset U$;
 
 \item $(AV'_I)$ for every $a\in A$ and every $U\in {\mathscr B}$ there exists $V\in{\mathscr B}$ such that $a.V\subset U$;
 
 \item $(AV'_{II})$ for every $U\in {\mathscr B}$ there exists $V\in {\mathscr B}$ such that $V.V\subset U$.
 
\end{description}
 \noindent
 Then, there exists a unique topology $\tau$ on $A$, compatible with the ring structure of $A$, such that ${\mathscr B}$ is a fundamental system of $\tau$-neighborhoods of $0$ .
 \end{proposition}
 
 \begin{proposition}\label {1.2} Let $A$ be a topological ring.
\begin{description}
 \item $(1^{\underline{\circ}})$ Let $X$ be an $A$-module and ${\mathscr B}$ a filter basis on $X$ verifying the conditions $(GA'_I)$ and $(GA'_{II})$ of Proposition \ref{1.1} and, in addition, the three conditions below:

\item $\,\,\,\,\,\,$$(MV'_I)$ for every $x_0\in X$ and every $V\in{\mathscr B}$ there exists a $0$-neighborhood $S$ in $A$ such that $S.x_0\subset V$;

\item $\,\,\,\,\,\,$ $(MV'_{II})$ for every $\lambda_0\in A$ and every $V\in{\mathscr B}$ there exists $W\in{\mathscr B}$ such that $\lambda_0.W\subset V$;

\item $\,\,\,\,\,\,$ $(MV'_{III})$ for every $V\in{\mathscr B}$ there exist $W\in{\mathscr B}$ and a $0$-neighborhood $S$ in $A$ such that $S.W\subset V$.

\end{description}
  \noindent
  Then, there exists a unique topology $\tau$ on $X$, compatible with the $A$-module structure of $X$, such that ${\mathscr B}$ is a fundamental system of $\tau$-neighborhoods of $0$.
\begin{description}
 \item $(2^{\underline{\circ}})$ Let $X$ an $A$-algebra and ${\mathscr B}$ a filter basis on $X$ verifying the conditions $(GA'_I)$, $(GA'_{II})$, $(MV'_I)$, $(MV'_{II})$ and $(MV'_{III})$ of $(1^{\underline{\circ}})$ and, in addition, the two conditions below:
 
 \item $\,\,\,\,\,\,$ $(AV'_I)$ for every $x_0\in X$ and every $V\in{\mathscr B}$ there exists $W\in{\mathscr B}$ such that $x_0.W\subset V$;

 \item $\,\,\,\,\,\,$ $(AV'_{II})$ for every $V\in{\mathscr B}$ there exist $W\in{\mathscr B}$ such that $W.W\subset V$.
 \end{description}
  \noindent
  Then, there exists a unique topology $\tau$ on $X$, compatible with the $A$-algebra structure of $X$, such that ${\mathscr B}$ is a fundamental system of $\tau$-neighborhoods of $0$.
 \end{proposition}

Let $X$ be a topological ring and $A$ a subring of $X$ (we assume here that the unit element of $X$ belongs to $A$, see \cite{L}). Then it is clear that the topology induced on $A$ by the topology on $X$ is compatible with the ring structure of $A$. Moreover, if we denote by $A_0$ the topological ring obtained by endowing $A$ with the topology induced by $X$, it is clear the $X$ becomes an $A_0$-topological algebra. Now, if ${\mathscr B}$ is a filter basis on $X$ verifying the conditions of Proposition \ref{1.1}, it is easily seen that
$$(AV'_I)\Rightarrow(MV'_I) \hskip0.3cm\hbox{and}\hskip0.3cm (MV'_{II}); \hskip0.5cm \hbox{and}\hskip0.5cm (AV'_{II})\Rightarrow (MV'_{III})$$
which allow the following simplification of Proposition \ref{1.2} ($2^{\underline{\circ}}$), which we shall need in the sequel:

\begin{corollary}\label{1.3} Let $X$ be a ring, $A$ a subring of $X$ {\rm(}see {\rm \cite{L})} and ${\mathscr B}$ a filter basis on $X$ verifying the four conditions $(GA'_I)$, $(GA'_{II})$, $(AV'_I)$ and $(AV'_{II})$ of Proposition \ref{1.1} {\rm (}with $X$ instead of $A${\rm )}. Then, there exists a unique topology $\tau$ on $X$ with the following properties:
\begin{description}
\item {\rm(I.)} $\tau$ is compatible with the ring structure of $X$. The topology $\tau_A$ induced by $\tau$ on $A$ is compatible with the ring structure of $A$;

\item {\rm(II.)} $X$ is an $A_0$-topological algebra {\rm(}here $A_0$ denotes the topological ring obtained by endowing $A$ with $\tau_A${\rm)}.

\item {\rm(III.)} ${\mathscr B}$ is a fundamental system of $\tau$-neighborhoods of $0$.
\end{description}
\end{corollary}


\section{The scalar full sharp topology}

 In the remainder of the paper we shall adopt  the following conventions. As usual, $\K$ denotes indistinctly $\real$ or $\comp$, $\I:=]0,1]$ and $\I_{\eta}:=]\,0,\eta[$ for each $\eta\in\I$. We shall assume fixed an arbitrary $m\in\nat^*$ and a non void open set $\Omega$ in $\real^m$. We shall use freely the symbols $A_q=A_q(m;\K)$ ($q\in\nat$) [A-B, Def.1.5 and Not.1.8], ${\mathscr G}(\Omega)$ [A-B, Not.2.1.1], $\OK=\OK(\real^m)$, ${\mathscr E}_M(\K)={\mathscr E}_M(\K;\real^m)$ and ${\mathscr N}(\K)={\mathscr N}(\K;\real^m)$ [A-B, 3.1, Def.3.1.2 and Rem.3.1.3]. The definitions of the algebra of the simplified generalized functions and the ring of the simplified generalized numbers, 
 $${\mathscr G}_s(\Omega):={{\mathscr E}_{M,s}[\Omega]\over{{\mathscr N}_s[\Omega]}}\,\,\,\hbox{ and }\,\,\,\OK_s:={{\mathscr E}_{M,s}(\K)\over{{\mathscr N}_s(\K)}}\,\,,$$
that we  adopt here are the ones given in [A-J, section 1], where the notation used are ${\mathscr G}(\Omega)$ and $\OK$. The only difference between  the above definitions in the simplified cases and the ones given in [A-B, 8.1 and 8.2] is  that here (resp. in [A-B]) the domains of $u\in{\mathscr E}_{M,s}[\Omega]$ and $v\in{\mathscr E}_{M,s}(\K)$ are $\I\times\Omega$ and $\I$ (resp. $\real^*_+\times\Omega$ and $\real^*_+$).

As usual, $K\subset\subset\Omega$ denotes the sentence ``$K$ is a compact subset of $\Omega$''.

We use an exhaustive sequence of open sets $(\Omega_l)_{l\in\nat}$ for $\Omega$, which means that $\Omega=\displaystyle\cup_{l\in\nat}\Omega_l$,   $\Omega_l$ being a non-void open set and $\overline{\Omega}_l\subset\subset\Omega_{l+1}$, for all $l\in\nat$. If $F:\Omega\longrightarrow\K$ is a function and $\emptyset\neq S\subset\Omega$, we define $||F||_S:=\displaystyle\sup_{x\in S}|F(x)|$ and, in the particular case $S=\Omega_l$ (for any $l\in\nat$) we write $||F||_l$ instead of $||F||_{\Omega_l}$, i.e., $||F||_l:=\sup_{x\in\Omega_l}|F(x)|$.

If $\varphi\in A_0=A_0(m,\K)$ we define $i(\varphi):=diam\,supp(\varphi)$ and since $0\notin A_0$ we have $i(\varphi)>0$ $\forall\varphi\in A_0$. It is easy to see that $i(\varphi_{\varepsilon})=\varepsilon i(\varphi)$ for all $\varepsilon >0$.

Our next result paves the way for the definition of a partial order relation on $\OR$ by using the same kind of ideas which leads to the definition of a partial order relation on $\,{\OR}_s$ in [A-J-O-S, sec.3], see also [G-K-O-S].

\begin{lemma}\label{2.1} For a given $x\in\OR$ the following are equivalent:
\begin{description}
\item
{\rm (i)} every representative $\widehat x$ of $x$ satisfies the condition:

\begin{math}
(*) \left |\begin{array}{l}
\exists N\in\nat \,\,\,\hbox{{\it such that}}\,\,\, \,\forall \,b>0\,\,\,\, \hbox{{\it and} } \,\,\,\forall\, \varphi\in A_N \,\,\,\,\hbox{\it{there is}}\\
\eta=\eta(b,\varphi)\in\I \,\,\,\, {\hbox{\it verifying }}\,\,\,\,\widehat x(\varphi_{\varepsilon})\geq -\varepsilon^{b}\,\,\,\, \forall \,\varepsilon\in\I_{\eta}\,;\end{array}\right.
\end{math}

\item {\rm (ii)} there exists a representative $\widehat x$ of $x$ satisfying $(*)$;

\item {\rm (iii)} there exists a representative $x_*$ of $x$ such that $x_*(\varphi)\geq 0\,\,\,\forall\,\varphi\in A_0$;

\item {\rm (iv)} there exist $N\in\nat$ and a representative $x_*$ of $x$ such that $x_*(\varphi)\geq 0\,\,\,\forall\,\varphi\in A_N$.

\end{description}
\end{lemma}

\begin{proof} $(ii)\Longrightarrow(iii)$: From the representative $\widehat x$ of $x$ in (ii) we define $h:A_0\longrightarrow\real$ by $h(\varphi):=0$ (resp. $-\widehat x(\varphi)$) if $\widehat x(\varphi)\geq 0$ (resp. $\widehat x(\varphi)<0$), hence $h(\varphi)\geq 0$ for each $\varphi\in A_0$ and it is easily seen that $h\in{\mathscr N}(\real)$. Now, it is enough to define $x_{*}(\varphi):=\widehat x(\varphi)+h(\varphi)\,\forall\,\varphi\in A_0$.

\noindent
$(iii)\Longrightarrow(i)$: Let $\widehat x$ be an arbitrary representative of $x$. Since $x_*-\widehat x\in{\mathscr N}(\real)$ it follows that there is $N\in\nat$ such that for all $q\geq N$, for all $\varphi\in A_q$ and for all $b>0$ we have ${|(x_*-\widehat x\,)(\varphi_{\varepsilon})|{\varepsilon^{-b}}}\rightarrow 0$ if $\varepsilon \downarrow 0$. Therefore, there is $\eta=\eta(b,\varphi)\in \I$ satisfying ${|(x_*-\widehat x\,)(\varphi_{\varepsilon})|{\varepsilon ^{-b}}}\leq 1\,$ $\forall\,\varepsilon\in \I_{\eta}$, hence $(x_*-\widehat x\,)(\varphi_{\varepsilon})\leq \varepsilon^b\,$ $\forall\,\varepsilon\in\I_{\eta}$, which implies $\widehat x(\varphi_{\varepsilon})\geq x_*(\varphi_{\varepsilon})-\varepsilon^b\geq -\varepsilon^b\,$ $\forall\,\varepsilon\in \I_{\eta}$, hence $\widehat x$ satisfies (*).

The implications $(i)\Longrightarrow(ii)$ and $(iii)\Leftrightarrow(iv)$ are obvious.\end{proof}

\begin{definition}\label{2.2} An element $x\in\OR$ is said to be quasi-positive or q-positive, if it has a representative satisfying the equivalent conditions of Lemma\ref{2.1}. We shall denote this by $x\geq 0$. We shall say also that $x$ is quasi-negative or q-negative if $-x$ is q-positive and we denote this by $x\leq 0$. If $y\in\OR$ is another element then we write $x\leq y \,($resp. $x\geq y )$ if $y-x$ $($resp. $x-y)$ is q-positive.
\end{definition}

\noindent
{\bf Notation.}

\vskip-0.5cm
$$\OK^{\,*}:=\OK\setminus\{0\}\,\,,\,\,\OR_+:=\{x\in\OR\,|\,x\geq 0\}\,\,\hbox{ and}\,\,\,\OR^{\,*}_+:=\{x\in\OR^{\,*}\,|\,x\geq 0\}.$$

\begin{example}\label{2.3} For every $r\in\real$ we shall define an element $\alpha_r^{\bullet}\in\OR^{\,*}_+$ which will play the same role that the elements $\alpha_r\in\OR_{s,+}$ $($see {\rm [A-J}, section 1$])$. A representative of $\alpha_r^{\bullet}$ is the function given by
$$\widehat{\alpha_r^{\bullet}}:\varphi\in A_0 \longmapsto i(\varphi)^r\in\real^*_+\,.$$
Clearly $\widehat{\alpha_r^{\bullet}}(\varphi)>0\,\forall\,\varphi\in A_0\,$ and $\widehat{\alpha_r^{\bullet}}(\varphi_{\varepsilon})=\varepsilon^r i(\varphi)^r=\widehat\alpha_r(\varepsilon)i(\varphi)^r$. Obviously $\widehat{\alpha_r^{\bullet}}\in{\mathscr E}_M(\real)$ and $\alpha_r^{\bullet}:=cl(\widehat{\alpha_r^{\bullet}})\in\OR_+^{\,*}$. From the condition $(*)$ of Lemma {\rm\ref{2.1}} it follows at once that

\vskip0.2cm
\begin{math}
(2.3.1)\begin{array}{l}
\,\,\,\,\,\,\,r,s\in\real\,\,\,\,\,\hbox{and}\,\,\,\,\,s<r \Rightarrow \alpha_r^{\bullet}\leq \alpha_s^{\bullet} \,.
\end{array}
\end{math}

\end{example}

\begin{proposition}\label{2.4} The relation $x\leq y$ of Definition {\rm\ref{2.2}} is a  partial order relation on $\OR$.
\end{proposition}
\begin{proof} Follows at once by remarking that from Lemma \ref{2.1}(i) we have that, for given $x\,,\,y\in\OR$, the relation $x\leq y$ means that
\vskip0.2cm
\begin{math}
(2.4.1)\left |\begin{array}{l}
\hbox{ for any representatives } \widehat x \hbox{ and } \widehat y \hbox{ of } x \hbox{ and } y,
\exists \,N\in\nat \\\hbox{ such that } \forall\, b>0 \hbox{ and } \forall\,\varphi\in A_N
\exists \,\eta=\eta(b,\varphi)\in\I \\\hbox{ verifying } \widehat x(\varphi_{\varepsilon})\leq \widehat y(\varphi_{\varepsilon})+\varepsilon^b\,\,\forall\,\varepsilon\in \I_{\eta}.
\end{array}\right.
\end{math}

\hskip11.5cm\end{proof}

\begin{example}\label{2.5} The partial order relation on $\OR$ it is not  linear. Indeed, by defining
$$\widehat x(\varphi):=i(\varphi)\sin\left ({1\over{i(\varphi)}}\right )\,\,\forall\, \varphi\in A_0$$
it is easily seen that $x:=cl(\widehat x)\geq 0$ and $x\leq 0$ are false.
\end{example}

\begin{proposition}\label{2.6} For every $x\in\OR_+$ and $p\in\nat^*$ there is  a unique $y\in\OR_+$ such that $y^p=x$ $($this $y$ is denoted by $x^{1/p}$ or $\root p\of x$ and it is called $q$-positive pth-root of $x$$)$.
\end{proposition}

\begin{proof}
Let $x_*$ be the representative of $x$ as in Lemma \ref{2.1}(iii). Then the function $y_*:\varphi\in A_0\longmapsto\root p\of {x_*(\varphi)}\in\real_+$ is well defined and moderate. Clearly $y:=cl(y_*)\in\OR_+$ and $y^p=x$. \end{proof}

For every $x\in\OK$, if $\widehat x$ is any representative of $x$, the function $|\widehat x|:\varphi\in A_0\longmapsto|\widehat x(\varphi)|\in\real_+$ is obviously moderate and since $|\widehat x|(\varphi)=|\widehat x(\varphi)|\geq 0\,\forall\,\varphi\in A_0$, it follows that $cl(|\widehat x|)\in\OR_+$. By the triangle inequality it follows that $cl(|\widehat x|)$ is independent of the representative $\widehat x$ of $x$ and only depends of $x$, hence it is natural to denote this class by $|x|$, i.e., 
$$|x|:=cl(|\widehat x|)$$
is called absolute value (or module) of $x$. So we have a natural map
$$x\in\OK\longmapsto |x|\in\OR_+.$$

The definition below depends on the above concept of absolute value and on Example \ref{2.3}.

\begin{definition}\label{2.7} If $x_0\in\OK$ and $r\in\real$, then we define
$$V_r[x_0]:=\{x\in\OK\,|\,|x-x_0|\leq\alpha_r^{\bullet}\}\,\,\,\,\,\,\,\,\hbox{and}\,\,\,\,\,\,\,\,{\mathscr B}:=\{V_r[0]\,|\,r\in\real\}.$$
\end{definition}

\begin{remark}\label{2.8} In view of Lemma {\rm\ref{2.1}(i)}, the statement $x\in V_r[0]$ is equivalent to the following sentence:

\vskip0.2cm
\begin{math}
(2.8.1) \left |\begin{array}{l}
\hbox{if}\,\, \widehat x \,\hbox{is any representative of }\,x\hbox{ then }\exists N\in\nat\,\hbox{ such that }\\ \,\forall\, b>0\,\hbox{ and }\,\,\,
\forall\,\varphi\in A_N\,\exists\,\eta=\eta(b,\varphi)\in\I\,\hbox{ verifying }\\
|\widehat x(\varphi_{\varepsilon})|\leq i(\varphi)^r\varepsilon^r+\varepsilon^b\,,\forall\,\varepsilon\in \I_{\eta}\,\end{array}\right.
\end{math}
\end{remark}

\begin{lemma}\label{2.9} The set ${\mathscr B}$ $($see Definition {\rm\ref{2.7}}$)$ is a filter basis on $\OK$ which satisfies the four conditions in Proposition {\rm\ref{1.1}}. Therefore, ${\mathscr B}$ determines a Hausdorff topology ${\mathscr T}$ compatible with the ring structure of $\OK$.
\end{lemma}
\begin{proof}
Clearly ${\mathscr B}\neq\emptyset$ and $\emptyset\notin{\mathscr B}$. The implication
$$r\,,\,s\in\real \Rightarrow \exists\,t\in\real\,\hbox{ such that }\, V_t[0]\subset V_r[0]\cap V_s[0]$$
follows trivially from the statement
$$r\,,\,s\in\real \,\hbox { and }\, r<s \Rightarrow  V_s[0]\subset V_r[0]\,,$$
which follows at once from (2.8.1).

Statement (2.8.1) implies also that ${\mathscr B}$ satisfies the four conditions of Proposition \ref{1.1}. Indeed, we get easily that $V_{r+1}[0]+V_{r+1}[0]\subset V_r[0]$ for all $r\in\real$, which shows that $(GA'_I)$ holds. The equality $V_r[0]=-V_r[0]\,\forall\,r\in\real$ shows that $(GA'_{II})$ holds. Now we must show that

\vskip0.2cm
\noindent
$(AV'_I)$ For given $a\in\OK$ and $V_r[0]\in{\mathscr B}$ there exists $s\in\real$ such that

\vskip0.2cm
\begin{math}
(2.9.1)\begin{array}{l}
\,\,\,\,\,\,\,\,\,\,\,\,\,\,\,\,\,\,\,\,\,\,\,\,\,\,\,a\,.\,V_s[0]\subset V_r[0]\,.\\
\end{array}
\end{math}

\vskip0.2cm
\noindent
In fact, since $a\in\OK$, if $\widehat a$ is any representative of $a$ $\,\exists\,N\in\nat\,$ such that for every $\varphi\in A_N\,\exists\,c>0\,$ and $\sigma\in \I$ satisfying $|\widehat a(\varphi_{\varepsilon})|\leq c\varepsilon^{-N}\,\forall\,\varepsilon\in \I_{\sigma}$. Then, it is easily seen that (2.9.1) holds by defining $s:=r+N+1\,$.

Next, we shall show that

\vskip0.2cm
\noindent
$(AV'_{II})$ Given $V_r[0]\in{\mathscr B}$ there is $s\in\real$ such that $V_s[0]\,V_s[0]\subset V_r[0]\,$, 

\vskip0.2cm
\noindent
which is trivial by defining $s:={r+1\over {2}}\,$.

Finally, $\mathscr T$ is Hausdorff since it is easily seen that $\,\displaystyle\cap_{r>0}V_r[0]=\{0\}$.
\end{proof}

\begin{definition}\label{2.10} The topology ${\mathscr T}$ on $\OK$ of Lemma {\rm\ref{2.9}} is called scalar full sharp topology.
\end{definition}

The definition of the partial order relation and the absolute value in $\OK$ just introduced in this section 2 are rather simple generalizations of the analogous definitions in $\OK_s$ developed in [A-J-O-S, Lemma 3.1, Definitions 3.2 and 3.4] that we shall use freely in the sequel.

Clearly, we can adapt Definition \ref{2.7} by introducing, for $x_0\in\OK_s$ and $r\in\real$ the set
$$V_r(x_0):=\{x\in\OK_s\,|\,|x-x_0|\leq \alpha_r\}\hskip0.5cm\hbox{and}\hskip0.5cm{\mathscr B}_s:=\{V_r(0)\,|\,r\in\real\}\,.$$

\begin{remark}\label{2.11} From {\rm [A-J-O-S, Lemma 3.1]} it is easily seen that the statement $x\in V_r(0)$ is equivalent to the following sentence:

\vskip0.2cm
\begin{math}
(2.11.1) \left |\begin{array}{l}
\hbox{if}\,\, \widehat x \,\,\hbox{is any representative of }\,\,x\,\hbox{ then } \,\forall\, b>0\,\,\exists\,\eta_b\in\I\,\\\hbox{ such }
\hbox{that }\,\,\,\,
 |\widehat x(\varepsilon)|\leq \varepsilon^r+\varepsilon^b\,,\,\forall\,\varepsilon\in \I_{\eta_b}\,\,.\hskip6.5cm\end{array}\right.
\end{math}
\end{remark}

With the above notation, introduced just before Remark \ref{2.11}, we have the following ``simplified version'' of  Lemma \ref{2.9}.

\begin{lemma}\label{2.12} The set ${\mathscr B}_s$ is a filter basis on $\OK_s$ which satisfies the four conditions of Proposition {\rm\ref{1.1}}. Therefore, ${\mathscr B}_s$ determine a topology ${\mathscr T}_s$ compatible with the ring structure of $\OK_s$.
\end{lemma}
\begin{proof}
The argument is a minor modification of the proof of Lemma 2.9 by using (2.11.1) instead (2.8.1).
\end{proof}

The topology ${\mathscr T}_s$ of Lemma 2.12 is indeed a familiar one:

\begin{theorem}\label{2.13} Let $\tau_s$ be the ``sharp topology'' on $\OK_s$ $($see {\rm \cite{S1}, \cite{S2}} and {\rm\cite{AJ}}$)$, then
\begin{description}
	\item {\rm(a)} ${\mathscr T}_s=\tau_s$;
	
	\item {\rm (b)} ${\mathscr T}_s$ coincides with the topology $j_m^{-1}({\mathscr T})$, where ${\mathscr T}$ was introduced in Definition {\rm\ref{2.10}} and 
	
	\vskip-0.6cm
	$$j_m:\lambda\in\OK_s\longmapsto cl(\widehat\lambda_*)\in\OK=\OK(\real^m)\,,$$
	being $\widehat\lambda$ an arbitrary representative of $\lambda$ and $\widehat\lambda_*(\varphi):=\widehat\lambda(min(1\,,\,i(\varphi)))\newline\forall\,\varphi\in A_0=A_0(\real^m\,;\,\K)\,$.
\end{description}
\end{theorem}
\noindent
Note that $j_m$ is the natural injective ring-homomorphism of $\OK_s$ into $\OK$, hence we can identifies $\OK_s$ with $Im(j_m)$ and write $\OK_s\subset\OK$. So we   can rewrite the statement {\rm (b)} of Theorem {\rm\ref{2.13}} by saying `` ${\mathscr T}_s=\tau_s$ is the topology induced by ${\mathscr T}$ ''.

\vskip0.2cm
\begin{proof}
(a) For given $x_0\in\OK_s$ and $\rho\in\real^*_+$ we set (see [A-J, Not.1.5])
$$B_{\rho}(x_0):=\{ x\in\OK_s\,|\,||x-x_0||<\rho \}\,.$$
Since the set of all balls $B_{\rho}(0)$ when $\rho$ ranges over $\real^*_+$ is a fundamental system of $\tau_s$-neighborhoods of $0$, it is enough to show the two following statements:

\vskip0.2cm
\begin{math}
(2.13.1)\begin{array}{l}
\,\,\,\,\,\,\forall\,r\in\real\,\exists\,\rho\in\real^*_+\,\hbox{ such that }\,B_{\rho}(0)\subset V_r(0)\\
\end{array}
\end{math}
\vskip0.2cm

\noindent
and

\vskip0.2cm
\begin{math}
(2.13.2)\begin{array}{l}
\,\,\,\,\,\,\forall\,\rho\in\real^*_+\,\exists\,r\in\real\,\hbox{ such that }\, V_r(0)\subset B_{\rho}(0)\,.\\
\end{array}
\end{math}
\vskip0.2cm

\noindent
In order to prove (2.13.1) fix an arbitrary $r\in\real$ and consider $\rho\in\real^*_+$ such that $\rho\leq e^{-r}$. We shall show that $B_{\rho}(0)\subset V_r(0)$. Indeed, 
$$x\in B_{\rho}(0)\Leftrightarrow ||x||=e^{-V(\widehat x)}<\rho\Leftrightarrow V(\widehat x)>-ln\,\rho\,,$$
where $\widehat x$ is any representative of $x$. The preceding inequality and the definition of $V(\widehat x)$ imply that there is $\eta\in \I$ such that
$$\displaystyle {|\widehat x(\varepsilon)|\over{\varepsilon^{-ln\,\rho}}}<1\,\,\forall\,\varepsilon\in \I_{\eta}\,.$$
Hence, from the choice of $\rho$, it follows that
$$|\widehat x(\varepsilon)|\leq \widehat\alpha_{-ln\,\rho}(\varepsilon)\leq\widehat\alpha_r(\varepsilon)\,\,\forall\,\varepsilon\in \I_{\eta}$$
which implies (see (2.11.1)) that $x\in V_r(0)\,$.

Let us now prove (2.13.2). To this end, fix an arbitrary $\rho\in\real^*_+$ and consider $r\in\real$ such that $e^{-r}<\rho$. We shall show that $V_r(0)\subset B_{\rho}(0)$ . In the proof of the preceding inclusion we shall need the following trivial statement:

\vskip0.2cm
\begin{math}
(2.13.3)\left |\begin{array}{l}
\hbox{ If } z\in\OK_s \hbox{ and } |z|\leq 1 \hbox{ then } V(\widehat z)\geq 0 \hbox{ for every}\\\hbox{ representative } \widehat z \hbox{ of } z\,.
\end{array}\right.
\end{math}

\vskip0.2cm
\noindent
Now, for a given $x\in\OK_s$ we have 
$$x\in V_r(0) \Leftrightarrow |x|\leq \alpha_r \Leftrightarrow |\alpha_{-r}x|\leq 1$$
and therefore, by (2.13.3), if $\widehat x$ is any representative of $x$ we get (see [A-J, Prop.1.3 (c)]) $V(\widehat{\alpha}_{-r}\,.\,\widehat x\,)=-r+V(\widehat x\,)\geq 0$ thus $V(\widehat x\,)\geq r$. Therefore, from the choice of $r$, we then conclude that $||x||=e^{-V(\widehat x\,)}\leq e^{-r}<\rho\,$, i.e., $x\in B_{\rho}(0)\,$.

(b) It is enough to show that $j_m^{-1}(V_r[0])=V_r(0)\,\forall\,r\in\real$, which follows, as usual, from (2.11.1) and (2.8.1).\end{proof}

\section{The full sharp topology}

We use again the notation of section 2. It will be convenient to give the following characterizations of the elements of ${\mathscr E}_M[\Omega]$ and ${\mathscr N}[\Omega]$ (of course, equivalent to the ones given in [A-B, Notation 2.1.1, (b) and (d)]):

\vskip0.3cm
\noindent
$u\in {\mathscr E}_M[\Omega]$ if and only if 

\vskip0.2cm
\begin{math}
[3.1] \left |\begin{array}{l}
\forall\,l\,,\,p\in\nat\,\,\exists\,N\in\nat\,\,\hbox{and}\,\,\exists\,\sigma\in\real\,
\,\hbox{ such that }\,\,\forall\,\varphi\in A_N\,\,\hbox{ we have}\\
||\partial^{\beta}u(\varphi_{\varepsilon}\,,\,.)||_l=o(\varepsilon^{\sigma})\,,\hbox{ if }\,\varepsilon\downarrow 0\,,\,\hbox{ whenever }\,\,|\beta|\leq p\,;
\end{array}\right.
\end{math}

\vskip0.5cm
\noindent
and $u\in {\mathscr N}[\Omega]$ if and only if 

\vskip0.2cm
\begin{math}
[3.2] \left |\begin{array}{l}
\forall\,l\,,\,p\in\nat\,\,\forall\,\sigma\in\real\,\,\exists\,N\in\nat\,\,\hbox{ such that }\,\,\forall\,\varphi\in A_N\,\,\,\,\hbox{we have}\\
||\partial^{\beta}u(\varphi_{\varepsilon}\,,\,.)||_l=o(\varepsilon^{\sigma})\,,\hbox{ if }\,\varepsilon\downarrow 0\,,\,\hbox{ whenever }\,\,|\beta|\leq p\,.
\end{array}\right.
\end{math}

\vskip0.5cm
We shall need the result below, whose trivial proof is omited.

\begin{lemma}\label{3.1} Let $F\,,\,G$ and $H$ be elements of ${\mathscr C}(\Omega;\K)$ such that 

\vskip-0.3cm
$$|\, |F(x)|-|G(x)|\,|\leq |H(x)|\,\,\,,\,\forall \,x\in\Omega\,.$$
Then, for each $K\subset\subset\Omega$, we have $\,|\,||F||_K-||G||_K\,|\leq ||H||_K\,$.
\end{lemma}

\begin{lemma}\label{3.2} For given $u\in{\mathscr E}_M[\Omega]\,,\,\beta\in\nat^m\,$ and $\,l\in\nat$, we have:
\begin{description}
\item{\rm (a)} the function
$$u^{\beta}_l:\varphi\in A_0\longmapsto ||\partial^{\beta}u(\varphi,.)||_l\in\real_+$$
is moderate {\rm (}i.e. $u^{\beta}_l\in{\mathscr E}_M(\real)${\rm )} and $cl(u^{\beta}_l)\in\OR_+\,$;

\item{\rm (b)} if $v\in{\mathscr E}_M[\Omega]$ and $u-v\in{\mathscr N}[\Omega]$, then $u^{\beta}_l-v^{\beta}_l\in{\mathscr N}(\real)$ and hence $cl(u^{\beta}_l)=cl(v^{\beta}_l)\,$.
\end{description}
\end{lemma}

\begin{proof}
The assertion (a) is clear and the statement (b) follows at once from definitions and Lemma \ref{3.1} .\end{proof}

The definition below makes sense from Lemma \ref{3.2}.

\begin{definition}\label{3.3} Fix $\beta\in\nat^m$ and $l\in\nat$. For every $f\in{\mathscr G}(\Omega)$ we define
$$||f||_{\beta,l}:=cl(\widehat f_l^{\,\beta})=cl[\varphi\in A_0\longmapsto ||\partial^{\beta}\widehat f(\varphi,.)||_l\in\real_+]$$
where $\widehat f$ is any representative of $f$. For every $f_0\in{\mathscr G}(\Omega)$ and $r\in\real$ we define
$$W^{\beta}_{l,r}[f_0]:=\{f\in{\mathscr G}(\Omega)\,|\,||f-f_0||_{\sigma,l}\leq\alpha_r^{\bullet},\,\hbox{ whenever }\,\sigma\leq\beta\}$$
and
$${\mathscr B}_{\Omega}:=\{W^{\beta}_{l,r}[0]\,|\,\beta\in\nat^m\,,\,l\in\nat\,\hbox{ and }\,\,r\in\real\}.$$
\end{definition}

For the sake of simplicity, in the statement of the result below, we consider that $\OK\subset{\mathscr G}(\Omega)$ as the generalized constants.

\begin{lemma}\label{3.4} For all $\beta\in\nat^m$, $l\in\nat$ and $r\in\real$ we have $\OK\cap W^{\beta}_{l,r}[0]=V_r[0]$ $($see Definition {\rm\ref{2.7}}$)$.
\end{lemma}

\begin{proof}
The inclusion $\OK\subset{\mathscr G}(\Omega)$ is given by the canonical homomorphism of $\OK$-algebras
$$\lambda=cl(\widehat\lambda)\in\OK\longmapsto \lambda^*:=cl[(\varphi,x)\in A_0\times\Omega \longmapsto \widehat{\lambda}(\varphi)\in\real_+]\in{\mathscr G}(\Omega)\,.$$
Therefore
$$\lambda\in\OK\cap W^{\beta}_{l,r}[0]\Leftrightarrow \lambda\in\OK\,\hbox{ and }\,||\lambda^*||_{\sigma,l}\leq\alpha^{\bullet}_r\,\,,\forall\,\sigma\leq\beta$$
and since from Definition \ref{3.3} we get (where $\delta_{ij}$ denotes the Kronecker $\delta$):
$$||\lambda^*||_{\sigma,l}=|\lambda|.\delta_{0,|\sigma|}\,(\forall\,\sigma\leq\beta)\,,$$
it follows that

\vskip0.2cm
\hskip1cm
$\lambda\in\OK\cap W^{\beta}_{l,r}[0]\Leftrightarrow \lambda\in\OK\,\hbox{ and }\,|\lambda|\leq\alpha^{\bullet}_r\Leftrightarrow \lambda\in V_r[0]\,.$\end{proof}

\begin{remark}\label{3.5} {\rm (a)} From Lemma {\rm\ref{2.1}(i)} it is easily seen that the statement $f\in W^{\beta}_{l,r}[0]$ is equivalent to the following sentence

\vskip0.2cm
\begin{math}
(3.5.1) \left |\begin{array}{l}
\hbox{for every representative}\,\widehat f \,\hbox{of }\,f, \exists\,N\in\nat\,\hbox{ such that }\,
\forall\, b>0\,\\\hbox{ and }\,
\forall\,\varphi\in A_N\,\exists\,\eta=\eta(b,\varphi)\in\I\,\hbox{ satisfying }\,\\
||\partial^{\sigma}\widehat f(\varphi_{\varepsilon},.)||_l\leq i(\varphi)^r\varepsilon^r+\varepsilon^b\,,\forall\,\varepsilon\in \I_{\eta}\,,\forall \,\sigma\leq\beta\,.\end{array}\right.
\end{math}

\vskip0.3cm
\noindent
{\rm (b)} The following three statements hold:
\begin{description}
\item {\rm (I.)} $\lambda\geq\lambda'\Rightarrow W^{\lambda}_{p,r}[0]\subset W^{\lambda'}_{p,r}[0]\,,\forall\, (p,r)\in\nat\times\real$ $($here, $\lambda\geq\lambda'$ means that $\lambda_j\geq\lambda'_j\,$, $\forall\,j=1,2,...,m)$.

\item {\rm (II.)} $r>r'\Rightarrow W^{\lambda}_{p,r}[0]\subset W^{\lambda}_{p,r'}[0]\,,\forall\,(\lambda,p)\in\nat^m\times\nat\,.$

\item{\rm (III.)} $p>q \Rightarrow W^{\lambda}_{p,r}[0]\subset W^{\lambda}_{q,r}[0]\,,\forall\,(\lambda,r)\in\nat^m\times\real\,.$
\end{description}
\end{remark}

\begin{theorem}\label{3.6} {\rm (a)} The set ${\mathscr B}_{\Omega}$ $($see Definition {\rm\ref{3.3}}$)$ is a filter basis on ${\mathscr G}(\Omega)$ which satisfies the four conditions $(GA'_I)\,,\,(GA'_{II})\,,\,(AV'_I)\,$ and $\,(AV'_{II})$ of Corollary {\rm\ref{1.3}}.

\vskip0.2cm
\noindent
{\rm (b)} There exists a unique topology ${\mathscr T}_{\Omega}$ on ${\mathscr G}(\Omega)$ verifying the following conditions:
\begin{description}
\item {\rm (I.)} ${\mathscr T}_{\Omega}$ is compatible with the ring structure of ${\mathscr G}(\Omega)$. The topology induced by ${\mathscr T}_{\Omega}$ on $\OK$ coincides with the topology ${\mathscr T}$ $($see Definition {\rm\ref{2.10}}$)$.

\item {\rm (II.)} ${\mathscr G}(\Omega)$ is a $\OK$-topological algebra $($when ${\mathscr G}(\Omega)$ and $\OK$ are endowed with ${\mathscr T}_{\Omega}$ and ${\mathscr T}$ respectively$)$.

\item {\rm (III.)} ${\mathscr B}_{\Omega}$ is a fundamental system of ${\mathscr T}_{\Omega}$-neighborhoods of $\,0\,$.

\item {\rm (IV.)} ${\mathscr T}_{\Omega}$ is metrizable.
\end{description}
\end{theorem}

\begin{proof} We abbreviate  notation by writing $W^{\beta}_{p,r}$ instead of $W^{\beta}_{p,r}[0]$. 

\noindent
(a) Clearly ${\mathscr B}_{\Omega}\neq\emptyset$ and $\emptyset\notin{\mathscr B}_{\Omega}$ since $0\in W^{\beta}_{p,r}$ for all $\beta\,,\,p\,,\,r\,.$ Given $W^{\beta}_{p,r}$ and $W^{\gamma}_{q,s}$ in ${\mathscr B}_{\Omega}$, if we define $\lambda=(\lambda_1,...,\lambda_m)$ where $\lambda_j:=\max(\beta_j,\gamma_j)\,(1\leq j\leq m)\,,\,l:=\max(p,q)\,$ and $\,t:=\max(r,s)\,$, from Remark \ref{3.5}(b), we get $W^{\lambda}_{l,t}\subset W^{\beta}_{p,r}\cap W^{\gamma}_{q,s}$. Hence, ${\mathscr B}_{\Omega}$ is a filter basis on ${\mathscr G}(\Omega)$. The verification of that ${\mathscr B}_{\Omega}$ satisfies the four conditions of Corollary \ref{1.3} is easy but rather tedious and here we only sketch these proofs.

\noindent
Condition $(GA'_I)$: Given any $W^{\beta}_{p,r}$ we have 
$$W^{\beta}_{p,s}+W^{\beta}_{p,s}\subset W^{\beta}_{p,r}\,\,\forall\,s>r\,.$$
Indeed, this follows easily from the implication
$$s>r\Rightarrow 2\alpha^{\bullet}_s\leq\alpha^{\bullet}_r\,,$$
which is a consequence of Lemma \ref{2.1}(i).

\noindent
Condition $(GA'_{II})$: Obvious, since $W^{\beta}_{p,r}=-W^{\beta}_{p,r}\,\,\forall\,\beta\,,\,p\,,\,r\,.$

\noindent
Condition $(AV'_I)$: Given $f_0\in{\mathscr G}(\Omega)$ and $W^{\beta}_{p,r}$ there exists $s\in\real$ such that

\noindent
(3.6.1)\hskip4cm $f_0W^{\beta}_{p,s}\subset W^{\beta}_{p,r}\,.$

\noindent
Fix a representative $\widehat f_0$ of $f$ then, from the moderation we can find $N_1\in\nat$ and $r'<0$ such that $\,\forall\varphi\in A_{N_1}\,\exists\,C>0$ and $\eta_1=\eta_1(\varphi)\in\I\,$ satisfying

\vskip0.2cm
\noindent
$(3.6.2)$\hskip2cm $||\partial^{\mu}\widehat f_0(\varphi_{\varepsilon},.)||_p\leq C\varepsilon^{r'}\,\,\forall\,\varepsilon\in\,I_{\eta_1}\,,\,\hbox{ whenever }\, |\mu|\leq |\beta|\,$

\vskip0.2cm
\noindent
Define $s:=r-r'+1=r+|r'|+1$. Now, the proof of (3.6.1) is easy: fix any representative $\widehat g$ of an arbitrary $g\in W^{\beta}_{p,s}$, we want to prove that $f_0g\in W^{\beta}_{p,r}$. From $g\in W^{\beta}_{p,s}$ and (3.5.1) we have an inequality for $||\partial^{\sigma}\widehat g(\varphi_{\varepsilon},.)||_p\,$ ($\sigma\leq\beta$) which jointly with (3.6.2) and the Leibnitz formula for derivation of a product, give an inequality for $||\partial^{\lambda}(\widehat f_0\widehat g)(\varphi_{\varepsilon},.)||_p\,$ ($\lambda\leq \beta)$, which shows (by (3.5.1)) that $f_0g\in W^{\beta}_{p,r}$.

\noindent
Condition $(AV'_{II})$: Given $W^{\beta}_{p,r}$ there is $s\in\real$ such that
$$W^{\beta}_{p,s}\,.\,W^{\beta}_{p,s}\subset W^{\beta}_{p,r}\,.$$
Define $s:=(r+1)/2\,.$ In order to prove the above inclusion, for fixed arbitrary $f\,,\,g\in W^{\beta}_{p,s}$ we want to proof that $fg\in W^{\beta}_{p,r}$. From $f\,,\,g\in W^{\beta}_{p,s}$ and (3.5.1), for fixed arbitrary representatives $\widehat f\,,\,\widehat g\,$ of $f\,,\,g$ respectively, we get two inequalities for $||\partial^{\sigma}\widehat f(\varphi_{\varepsilon},.)||_p\,$ and $||\partial^{\sigma}\widehat g(\varphi_{\varepsilon},.)||_p\,$ ($\sigma\leq\beta)$. Next, we apply Leibnitz formula which, jointly with the above inequalities, leads to an inequality for $||\partial^{\lambda}(\widehat f\,\widehat g)(\varphi_{\varepsilon},.)||_p\,$ ($\lambda\leq\beta$), which shows (by (3.5.1)) that $f\,g\in W^{\beta}_{p,r}\,$.

\noindent
(b) In view of (a) and Corollary \ref{1.3} we can conclude that there exists a unique topology ${\mathscr T}_{\Omega}$ on ${\mathscr G}(\Omega)$ satisfying the three following conditions:

(I'.) ${\mathscr T}_{\Omega}$ is compatible with the ring structure of ${\mathscr G}(\Omega)$. The topology ${\mathscr T}_{\OK}$ induced by ${\mathscr T}_{\Omega}$ on $\OK$, is compatible with the ring structure of $\OK$;

(II') ${\mathscr G}(\Omega)$ is a $\OK$-topological algebra $({\mathscr G}(\Omega)$ and $\OK$ are endowed with ${\mathscr T}_{\Omega}$ and ${{\mathscr T}}_{\OK}$ respectively);

(III') ${\mathscr B}_{\Omega}$ is a fundamental system of ${\mathscr T}_{\Omega}$-neighborhoods of $0$.

\noindent
Now, by Lemma \ref{3.4}, it is clear that the topology induced by ${\mathscr T}_{\Omega}$ on $\OK$ is ${\mathscr T}$ (see Definition \ref{2.10}), hence ${\mathscr T}_{\OK}={\mathscr T}$ , and therefore (I'.) implies (I.) and (II'.) implies (II.).

Since $\lim\limits_{r\to+\infty}\alpha_r^{\bullet}=0\,$ it is clear that the set

\vskip-0.4cm
 $$\{W^{\beta}_{p,r}\,|\,\beta\in\nat^m\,,\,p\in\nat\,\hbox{ and }\,r\in\nat\}$$
 
 \vskip-0.2cm
 \noindent
is a countable fundamental system of ${\mathscr T}_{\Omega}$-neighborhoods of $\,0\,$ and so (IV) is true.\end{proof}

\begin{definition}\label{3.7} The topology ${\mathscr T}_{\Omega}$ on ${\mathscr G}(\Omega)$ in Theorem {\rm\ref{3.6}} is called full sharp topology.
\end{definition}

In the sequel, we shall show that the topology induced by ${\mathscr T}_{\Omega}$ on ${\mathscr G}_s(\Omega)$ (see notation in section 2) is precisely the sharp topology $\tau_{\Omega}$  (see \cite{S1}, \cite{S2} and \cite{AJ}). Here, we shall develop a procedure analogous to the one at the beginning  of section 3, proving the analog  of Lemma \ref{3.2} for ${\mathscr G}_s(\Omega)$, which allows us to define a filter basis ${\mathscr B}_{s,\Omega}$ on ${\mathscr G}_s(\Omega)$ which will determine a topology ${\mathscr T}_{s,\Omega}$ on ${\mathscr G}_s(\Omega)$. Next, we shall show that
$$W^{\beta}_{p,r}[0]\cap {\mathscr G}_s(\Omega)=W^{\beta}_{p,r}(0)\,\,\,\,\forall\,\beta\,,\,p\,,\,r$$
where $W^{\beta}_{p,r}(0)$ denote a general element of ${\mathscr B}_{s,\Omega}$. Hence ${\mathscr T}_{s,\Omega}$ will be the topology induced on ${\mathscr G}_s(\Omega)$ by ${\mathscr T}_{\Omega}$. Finally, we shall show that ${\mathscr T}_{s,\Omega}=\tau_{\Omega}$.

\begin{lemma}\label{3.8} Given $u\in{\mathscr E}_{M,s}[\Omega]$, $\beta\in\nat^m$ and $l\in\nat$, we have
\begin{description}
\item {\rm (a)} the function

\vskip-0.7cm
$$u^{\beta}_l:\varepsilon\in\I\longmapsto ||\partial^{\beta}u(\varepsilon,.)||_l\in\real_+$$
is moderate $($i.e. $u^{\beta}_l\in{\mathscr E}_{M,s}(\real))$ and $cl(u^{\beta}_l)\in\OR_{s,+}\,$;

\item{\rm (b)} if $v\in{\mathscr E}_{M,s}[\Omega]$ and $u-v\in{\mathscr N}_s[\Omega]$, then $u_l^{\beta}-v_l^{\beta}\in {\mathscr N}_s(\real)$ and therefore $cl(u^{\beta}_l)=cl(v^{\beta}_l)$.
\end{description}
\end{lemma}\begin{proof} Follows at once from slight modifications in the proof of Lemma \ref{3.2}.\end{proof}

\begin{definition}\label{3.9} Fix $\beta\in\nat^m$ and $l\in\nat$. For every $f\in{\mathscr G}_s(\Omega)$ we define
$$|f|_{\beta,l}=cl({f_*}^{\beta}_l)=cl[\varepsilon\in \I\longmapsto ||\partial^{\beta}f_*(\varepsilon,.)||_l\in\real_+]\,,$$
where $f_*$ is any representative of $f$. For every $f_0\in{\mathscr G}_s(\Omega)$ and $r\in\real$ we define
$$W^{\beta}_{l,r}(f_0):=\{f\in{\mathscr G}_s(\Omega)\,|\,|f-f_0|_{\sigma,l}\leq \alpha_r\,\hbox{ whenever }\, \sigma\leq\beta\}$$
and
$${\mathscr B}_{s,\Omega}:=\{W^{\beta}_{l,r}(0)\,|\,\beta\in\nat^m\,,\,l\in\nat\,\hbox{ and }\, r\in\real\}\,.$$
\end{definition}

Here we have also  remarks analogous to Remark \ref{3.5}.

\begin{remark}\label{3.10} {\rm (a)} From {\rm [A-J-O-S}, Lemma {\rm 3.1(i)]}  it follows at once that the statement $f\in W^{\beta}_{l,r}(0)$ is equivalent to the following sentence 

\vskip0.2cm
\begin{math}
(3.10.1) \left |\begin{array}{l}
\hbox{every representative}\, f_* \,\hbox{of }\,f\, { satisfies}:
\,\forall\, b>0\,\hbox{ there exists }\\
\eta_b\in\I\,\hbox{ such that}\,\,
||\partial^{\sigma} f_*({\varepsilon},.)||_l\leq \varepsilon^r+\varepsilon^b\,,\forall\,\varepsilon\in \I_{\eta_b},\forall \,\sigma\leq\beta.\end{array}\right.
\end{math}

\vskip0.3cm
\noindent
{\rm (b)} The three implications of Remark {\rm\ref{3.5}(b)} hold for the sets $W^{\beta}_{l,r}(0)$, we omit them since it is suffices to change $W^{\beta}_{l,r}[0]$ by  $W^{\beta}_{l,r}(0)$.
\end{remark}

Since ${\mathscr G}_s(\Omega)\subset{\mathscr G}(\Omega)$, for given $f\in{\mathscr G}_s(\Omega)$, $\beta\in\nat^m$ and $l\in\nat$, the two generalized numbers $|f|_{\beta,l}\in\OK_s$ and $||f||_{\beta,l}\in\OK$ make sense. Moreover, since $\OK_s\subset\OK$, it is natural to ask the relation  between these two generalized numbers. We shall show a nice answer to this question

\vskip0.2cm
\begin{math}
\begin{array}{l}

[3.3] \hskip2cm |f|_{\beta,l}=||f||_{\beta,l}\,\,\,\,\,(\forall\,f\,,\,\beta\,,\,l)\,.
\end{array}
\end{math}

\vskip0.2cm
More precisely, if $j_m:\OK_s\longrightarrow \OK$ denote the canonical inclusion of $\OK_s$ into $\OK$ (see Theorem \ref{2.13}(b)) and 
$$\Psi^*_{\Omega}:f\in{\mathscr G}_s(\Omega)\longmapsto cl[\Psi_{\Omega}(f_*)]\in{\mathscr G}(\Omega)$$
denotes the canonical inclusion of ${\mathscr G}_s(\Omega)$ into ${\mathscr G}(\Omega)$ [where $f_*\in{\mathscr E}_{M,s}[\Omega]$ is any representative of $f$ and $\Psi_{\Omega}:u\in{\mathscr E}_{M,s}[\Omega]\longmapsto u_1\in{\mathscr E}_M[\Omega]$ is the ring homomorphism defined by $u_1(\varphi,x):=u(\min(1,i(\varphi)),x)\,$ for each $(\varphi,x)\in A_0\times\Omega$] then we have:

\begin{lemma}\label{3.11} {\rm (a)} $j_m(|f|_{\beta,l})=||cl[\Psi_{\Omega}(f_*)]||_{\beta,l}$ for every $f\in{\mathscr G}_s(\Omega)$, $\beta\in\nat^m$, $l\in\nat$ and any representative $f_*$ of $f$.

\vskip0.2cm
\noindent
{\rm (b)} $|f|_{\beta,l}\leq \alpha_r$ if and only if $||cl[\Psi_{\Omega}(f_*)]||_{\beta,l}\leq\alpha^{\bullet}_r$, for every $f\in{\mathscr G}_s(\Omega)$, $\beta\in\nat^m$, $l\in\nat$ and a representative $f_*$ of $f$.

\vskip0.2cm
\noindent
{\rm (c)} $\Psi{^*_{\Omega}}^{-1}(W^{\beta}_{p,r}[0])=W^{\beta}_{p,r}(0)$ for every $\beta\in\nat^m$, $p\in\nat$ and $r\in\real$.
\end{lemma}

\begin{proof} (a) Fix $f,\beta,l$ and any representative $f_*\in{\mathscr E}_{M,s}[\Omega]$ of $f$. The equality follows at once from the definitions of $\,|f|_{\beta,l}\,,\,\,j_m\,,\,\Psi_{\Omega}\,$ and  $||cl[\Psi_{\Omega}(f_*)]||_{\beta, l}$.
\newline
\noindent
(b) If $f_*$ is a representative of $f$, the inequality $|f|_{\beta,l}\leq\alpha_r$ is equivalent to the statement (3.10.1). Analogously, $||cl[\Psi_{\Omega}(f_*)]||_{\beta,l}\leq\alpha_r^{\bullet}$ is equivalent to the statement (3.5.1) for $\widehat f:=\Psi_{\Omega}(f_*)$. Now, the equivalence between the two inequalities follows as a tedious application of (3.10.1) and (3.5.1).

\noindent
(c) Follows immediately from (b).\end{proof}

Note that {\rm [3.3]} is an abuse of notation whose correct meaning is given in Lemma {\rm\ref{3.11}(a)}. Analogously, the relation

\vskip0.3cm
\noindent
{\rm [3.4]} \hskip3cm $W^{\beta}_{p,r}[0]\cap{\mathscr G}_s(\Omega)=W^{\beta}_{p,r}(0)\,\,\,\,(\forall\, \beta\,,\,p\,,\,r)$

\vskip0.3cm
\noindent
is an abuse of notation whose correct meaning is given in Lemma {\rm\ref{3.11}(c)}.

\begin{theorem}\label{3.12} {\rm (a)} The set ${\mathscr B}_{s,\Omega}$ $($see Definition {\rm\ref{3.9}}$)$ is a filter basis on ${\mathscr G}_s(\Omega)$ which satisfies the four axioms $(GA'_I)$, $(GA'_{II})$, $(AV'_I)$ and $(AV'_{II})$ of Corollary {\rm \ref{1.3}}.

\vskip0.2cm
\noindent
{\rm (b)} There exist a unique topology ${\mathscr T}_{s,\Omega}$ on ${\mathscr G}_s(\Omega)$ verifying the following conditions
\begin{description}
 \item{\rm (I.)} ${\mathscr T}_{s,\Omega}$ is compatible with the ring structure of ${\mathscr G}_s(\Omega)$. The topology induced by ${\mathscr T}_{s,\Omega}$ on $\OK_s$ coincides with the topology ${\mathscr T}_s$ $($see Lemma {\rm\ref{2.12}}$)$ and hence, with the sharp topology $\tau_s$ $($see Theorem {\rm\ref{2.13}}$)$;
 
\item{\rm (II.)} ${\mathscr G}_s(\Omega)$ is a $\OK_s$-topological algebra $($here we assume that ${\mathscr G}_s(\Omega)$ and $\OK_s$ are endowed with the topologies ${\mathscr T}_{s,\Omega}$ and ${\mathscr T}_s$ respectively$)$;

\item{\rm (III.)} ${\mathscr B}_{s,\Omega}$ is a fundamental system of ${\mathscr T}_{s,\Omega}$-neighborhoods of $0$.

\item{\rm (IV.)} ${\mathscr T}_{s,\Omega}$ coincides with the sharp topology $\tau_{\Omega}$;

\item{\rm (V.)} ${\mathscr T}_{s,\Omega}$ coincides with the topology induced by ${\mathscr T}_{\Omega}$ on ${\mathscr G}_s(\Omega)$.
\end{description}
\end{theorem}
\begin{proof} By identifying ${\mathscr G}_s(\Omega)$ with his image ${\mathscr G}(\Omega)$ by the canonical map $\Psi_{\Omega}^*$ we can write ${\mathscr G}_s(\Omega)\subset{\mathscr G}(\Omega)$ and then, the Lemma \ref{3.11} (c) shows that

\vskip0.2cm
\begin{math}
(3.12.1)\begin{array}{l}
\hskip1.5cm W^{\beta}_{p,r}[0]\cap{\mathscr G}_s(\Omega)=W^{\beta}_{p,r}(0)\,\,\,\forall\,\, \beta,p\,,\,r
\end{array}
\end{math}

\vskip0.2cm
\noindent
which means that

\vskip0.2cm
\begin{math}
(3.12.1')\begin{array}{l}
\hskip1.5cm {\mathscr B}_{s,\Omega}=\{W\cap {\mathscr G}_s(\Omega)\,|\,W\in {\mathscr B}_{\Omega}\,\}\,.
\end{array}
\end{math}

\vskip0.2cm
\noindent
The above relation (3.12.1'), together with Theorem \ref{3.6} (a), implies (a). Now, from (a) and Corollary \ref{1.3} we get (b), the first statement of (I.), (II.) and (III.). The proof of the second statement of (I.) follows at once by noting that
$$W^{\beta}_{p,r}(0)\cap\OK_s=V_r(0)\,\,\,\forall\,\,\beta,p\,,\,r$$
which is the simplified version of Lemma \ref{3.4} and easily proved. Note also that (V.) follows directly from (3.12.1').
\newline
(IV): Let recall that if $u\in{\mathscr E}_{M,s}[\Omega]$ and $(n,p)\in\nat^2$, we set

\vskip-0.2cm
 $$S_{np}(u):=\{\,a\in\real\,|\,||\partial^{\beta}u(\varepsilon\,,\,.)||_n=o(\varepsilon^a)\,\,\hbox{ if }\,\,\varepsilon\downarrow 0\,\,,\,\,\forall\,\,|\beta|\leq p\,\}$$
 
\noindent
and $v_{np}:=\sup\,S_{np}(u)$. Moreover, it is easy to see that $v_{np}$ is constant on every equivalence class $u+{\mathscr N}_s[\Omega]$, hence for each $f\in{\mathscr G}_s(\Omega)$  the (extended) real number

\vskip-0.4cm
$$V_{np}(f):=v_{np}(f_*)\in \,]-\infty,\infty\,]$$

\vskip-0.1cm
\noindent
is well defined, where $f_*$ is any representative of $f$. The sharp topology $\tau_{\Omega}$ on ${\mathscr G}_s(\Omega)$ is defined by the family of pseudo metrics

\vskip-0.4cm
$$d_{np}(f,g):=\exp(-V_{np}(f-g))\,\,\,\forall\,\,f,g\in{\mathscr G}_s(\Omega)\,,\,\,\forall\,\,n,p\in\nat\,.$$

\vskip-0.1cm
\noindent
For each $a>0$ we define the $d_{np}$-ball of center $0$ and radius $a$ :

\vskip-0.4cm
$$\overline{B}_a(d_{np}):=\{\,f\in{\mathscr G}_s(\Omega)\,|\,d_{np}(f,0)\leq a\,\}\,.$$

\vskip-0.1cm
\noindent
Then, the collection of all finite intersections of these balls is a fundamental system of $\tau_{\Omega}$-neighborhoods of $\,0$ in ${\mathscr G}_s(\Omega)$. Therefore it suffices to prove the two following statements:

\vskip0.2cm
\begin{math}
(3.12.2)\left |\begin{array}{l}
\hbox{ For every } a\in\real^*_+ \hbox{ and } (n,p)\in\real^2 \hbox{ there is a finite }\\ \hbox{ sequence } (W_{l,r}^{\beta}(0))_{\beta\in B}\,\, \hbox{ of elements in } \,\,{\mathscr B}_{s,\Omega} \hbox { such that }\\ \,\,\cap_{\beta\in \,B}W^{\beta}_{l,r}(0)\subset\overline {B}_a(d_{np})\,,\,.
\end{array}\right.
\end{math}

\vskip0.2cm
\noindent
and

\vskip0.2cm
\begin{math}
(3.12.3)\left |\begin{array}{l}
\hbox{ For each } \beta\in\nat^m, l\in\nat \hbox{ and } r\in\real  \hbox{ there is } \overline B_a(d_{np}) \\ \hbox { such that } \,\,\overline B_a(d_{np})\subset W^{\beta}_{l,r}(0).
\end{array}\right.
\end{math}

\vskip0.2cm

To obtain (3.12.2) fix $a\in\real^*_+$ and $(n,p)\in\nat^2$ as in (3.12.2).  Then, clearly we can supose without lost of generality that $0<a<e^{-1}$. Now, by defining $l:=n$ and choosing $r\in\real$ such that $r\geq -log\,a+1$, if $B:=\{\beta\in\nat^m\,|\,|\beta|=p\,\}$ it follows (apply (3.10.1) as usual) that $\cap_{\beta\in\, B}W^{\beta}_{n,r}(0)\subset\overline {B}_a(d_{np})\,$.

To prove (3.12.3) fix $\beta\in\nat^m$, $l\in\nat$ and $r\in\real$ as in (3.12.3). If we define $\,\,n:=l,p:=|\beta|$ and  take $a\in\,]\,0,e^{-r}\,[$, one proves easily the inclusion $\overline B_a(d_{lp})\subset W^{\beta}_{l,r}(0)$.\end{proof}


\section{A few results of convergence}

\begin{proposition}\label{5.1}$(a)$ Let $M$ be a Lebesgue-measurable set such that $\overline M\subset\subset\Omega$. Then the $\OK$-linear function
$$J_M: f\in{\mathscr G}(\Omega)\longmapsto \displaystyle\int_M\,f\,\in\OK$$
is ${\mathscr T}_{\Omega}-{\mathscr T}-$continuous;

\noindent
$(b)$ For every $\alpha\in\nat^m$ the $\OK$-linear function
$$\partial^{\alpha}:f\in{\mathscr G}(\Omega)\longmapsto \partial^{\alpha}f\in{\mathscr G}(\Omega)$$
is ${\mathscr T}_{\Omega}-{\mathscr T}_{\Omega}-$continuous.
\end{proposition}

\begin{proof}
(a) Fix $V_r[0]$ with $r\in\nat$ arbitrary. Then, for $\beta:=0=(0,0,...,0)\in\nat^m\,,\, s:=r+1$ and $l\in\nat$ such that $\overline M\subset\Omega_l$ one proves easily that
$$J_M(W_{l,r+1}^0[0])\subset V_r[0]\,.$$
(b) For  $W_{l,r}^{\beta}[0]$ given arbitrarily it is clear that $\partial^{\alpha}(W_{l,r}^{\alpha+\beta}[0])\subset W_{l,r}^{\beta}[0]\,.$
\end{proof}

\begin{corollary}\label{5.2} If $\,\,P=\displaystyle\sum_{|\alpha|\leq m} a_{\alpha}(x)\partial^{\alpha}\,\,$ is a  generalized LPDO $($i.e. $a_{\alpha}\in{\mathscr G}(\Omega)\,\,$ $\forall\,|\alpha|\leq m\,)$, then $P$ defines a linear application which is ${\mathscr T}_{\Omega}-{\mathscr T}_{\Omega}$- continuous
$$P:f\in{\mathscr G}(\Omega) \longmapsto Pf\in{\mathscr G}(\Omega)\,.$$
\end{corollary}

\begin{proof} If is enough to show that if $(f_l)_{l\in\nat}$ is a sequence in ${\mathscr G}(\Omega)$ such that  $f_l \stackrel{{\mathscr T}_{\Omega}}{\longrightarrow} f\in{\mathscr G}(\Omega)$ then $Pf_l\stackrel{{\mathscr T}_{\Omega}}{\longrightarrow}Pf$.

The continuity of the multiplication in ${\mathscr G}(\Omega)$ proves that for every $a\in{\mathscr G}(\Omega)$ the function $M_a:f\in{\mathscr G}(\Omega)\longmapsto af\in{\mathscr G}(\Omega)$ is continuous. From Proposition \ref{5.1} (b) we get $\,\partial^{\alpha}f_l \stackrel{{\mathscr T}_{\Omega}}{\longrightarrow} \partial^{\alpha}f\,$ and hence $M_a(\partial^{\alpha}f_l) \stackrel{{\mathscr T}_{\Omega}}{\longrightarrow} M_a(\partial^{\alpha}f)\,$ which implies that $\,\,a_{\alpha}\partial^{\alpha}f_l \stackrel{{\mathscr T}_{\Omega}}{\longrightarrow} a_{\alpha}\partial^{\alpha}f\,\,\,\forall\,\,|\alpha|\leq m\,.$ Now, the continuity of the  addition in ${\mathscr G}(\Omega)$ proves that $Pf_l \stackrel{{\mathscr T}_{\Omega}}{\longrightarrow} Pf\,$.\end{proof}

\begin{corollary}\label{5.3} If $\Omega$ is an open subset of $\,\,\comp^m$ then the subalgebra ${\mathscr H}{\mathscr G}(\Omega)\,$ of all holomorphic generalized functions on $\Omega$ is  a  ${\mathscr T}_{\Omega}$-closed subalgebra of $\,{\mathscr G}(\Omega)$.
\end{corollary}
\begin{proof} Let $(f_l)$ be a sequence in ${\mathscr H}{\mathscr G}(\Omega)$ and assume that $f_l\stackrel{{\mathscr T}_{\Omega}}{\longrightarrow}f\in{\mathscr G}(\Omega)$. From  Corollary \ref{5.2} above we have, for each $j=1,2,...,m$:

\vskip-0.2cm
$$0={\partial f_l\over{\partial\overline z_j}}\longrightarrow{\partial f\over{\partial\overline z_j}}\,\,,\,\,\hbox{ if }\,\,\,l\rightarrow \infty\,$$

\vskip-0.2cm
\noindent
hence $\,\displaystyle{\partial f\over{\partial\overline z_j}}=0$ for every $j=1,2,...,m$. \end{proof}


\section*{Appendix: On $\OK$-locally convex modules ($\OK$-LCM)}

With the notation introduced in Definition \ref{3.3} it is clear that for given $\beta\in\nat^m$ and $l\in\nat$ we have
$$||f+g||_{\beta,l}\leq ||f||_{\beta,l}+||g||_{\beta,l}\,\,\forall\,f\,,\,g\in{\mathscr G}(\Omega);$$
$$||a\,f||_{\beta,l}=|a|\,||f||_{\beta,l}\,\,\forall\,f\in{\mathscr G}(\Omega)\,\,\hbox{ and }\,\,a\in\OK\,,$$
which suggest that we can try to mimic some basic facts of the general theory of LCS by defining:

\vskip0.5cm
\noindent
{\bf Definition A.1}\label{A1} {\it Let $E$ be a ${\OK}$-module. A generalized semi-norm $($or a $G$-seminorm, for short$)$ on $E$ is a function $p:E\longrightarrow \OR_+$ verifying the two conditions:
\begin{description}
\item $(GSN1)$ $p(x+y)\leq p(x)+p(y)\,\,\,\forall\,x\,,\,y\in E$;

\item $(GSN2)$ $p(a\,x)=|a|\,p(x)\,\,\,\forall\,x\in E\,,\,\forall\,a\in\OK\,.$
\end{description}}

\vskip0.5cm
With the above notation clearly we have

\vskip-0.4cm
$$|p(x)-p(y)|\leq p(x-y)\,\,\,\forall\,x\,,\,y\in E$$

\vskip-0.1cm
\noindent
since the classical proof works in this case.

The remark preceding Definition A.1 shows that for every $\beta\in\nat^m$ and $l\in\nat$, the function

\vskip-0.5cm
$$||\,.\,||_{\beta,l}:f\in{\mathscr G}(\Omega)\longmapsto ||f||_{\beta,l}\in\OR_+$$

\vskip-0.2cm
\noindent
is a $G$-seminorm.

Next, fix a $\OK$-module $E$ and a $G$-seminorm $p$ on $E$. For each $r\in\real$ we can define the $p$-ball of center $0$ and radius $\alpha_r^{\bullet}$ by

\vskip-0.3cm
$$B_{p,r}=B_{p,r}(0):=\{x\in E\,|\,p(x)\leq\alpha_r^{\bullet} \}$$

\noindent
and for $x_0\in E$ we set $B_{p,r}(x_0):=x_0+B_{p,r}=\{x\in E\,|\,p(x-x_0)\leq\alpha_r^{\bullet}\}$. We can also  define the generalized segment of extremities $0$ and $1$ by

\vskip-0.3cm
$$[0,1]_g:=\{\lambda\in\OR\,|\,0\leq\lambda\leq 1\}\,,$$

\vskip-0.1cm
\noindent
note that $\alpha^{\bullet}_r\in [0,1]_g$ for every $r\in\R_+$.

A subset $A$ of $($the $\OK$-module$)$ $E$ is $G$-convex $($resp. $G$-absolutely convex$)$ if $x\,,\,y\in A$, then $\lambda x+(1-\lambda)y\in A$ for all $\lambda\in [0,1]_g$ $($resp. if $x\,,\,y\in A$, then $\alpha x+\beta y\in A$ for all $\alpha\,,\,\beta\in \OK$ such that $|\alpha|+|\beta|\leq 1 )$. So the definition below makes sense.

\vskip0.5cm
\noindent
{\bf Definition A.2}\label{A.2} {\it Let $E$ be a $\OK$-module. A module topology ${\mathscr T}$ on $E$ $($i.e. a topology compatible with $\OK$-module structure$)$ is said to be $G$-locally convex if $\,0$ has a fundamental system of $G$-absolutely convex neighborhoods.}

\vskip0.5cm
Let $A$ be a $\K$-algebra and $\Gamma$ a non void set of seminorms in $A$ such that: (I.) $p(x)=0\,\,\forall\,p\in\Gamma \Leftrightarrow x=0\,$; (II.) For all $p_1\,,\,p_2\in\Gamma$ there is $q\in\Gamma$ such that $\,p_i\leq q\,(i=1\,,\,2)$; (III.) For each $p\in\Gamma$ and each $\alpha>0$ we have $\alpha p\in\Gamma$. By defining $B_p:=\{x\in A\,|\,p(x)\leq 1\,\}\,\, \forall\,\,p\in\Gamma$, one proves trivially that the continuity of the multiplication 

\vskip-0.4cm
$$(x,y)\in A\times A\longrightarrow xy\in A$$

\vskip-0.2cm
\noindent
is  equivalent to each of the two following conditions:

\vskip0.2cm
 (i) $\forall p\in\Gamma\,\,\exists\,q\in\Gamma\,\,\hbox{ such that } \,\,B_qB_q\subset B_p$ ;

\vskip0.2cm
(ii) $\forall p\in\Gamma\,\,\exists\,q\in\Gamma\,\,\hbox{ such that } \,\, p(xy)\leq q(x)q(y)\,\,\forall\,x\,,\,y\in A$.

\vskip0.2cm
\noindent
In this case, the Hausdorff topological algebra $(A,\Gamma)$ is said to be a {\it locally multiplicatively-convex algebra} (see \cite{H}, in this book the condition (ii) has a slight mistake). Note that the topology ${\mathscr T}_{\Gamma}$ determinated by $\Gamma$ on $A$ is given (see Corol.\ref{1.3}) by the filter basis ${\mathscr B}_{\Gamma}$ of all finite intersection of balls $B_p$.

In our case $(\OK$ instead of $\K)$, if $A$ is a $\OK$-algebra and $\Gamma$ is a non-void set of $G$-seminorms on $A$ it is easy to prove that, for the conditions (i), (ii) above, we have $\hbox{(ii)} \Longrightarrow \hbox {(i)}$ but, seemingly $(\hbox{i}) {\Longrightarrow\kern -1.5 em /} \,\,\,\,\, (\hbox{ii})$. Therefore, it is natural, in our case, to say that the topological algebra $(A,\Gamma)$ is a {\it $G$-locally multiplicatively-convex algebra} if $\Gamma$ satisfies the above condition (i). Note that the proof of Theorem \ref{3.6} (a) $(AV'_{II})$ shows that $({\mathscr G}(\Omega),\Gamma)=({\mathscr G}(\Omega),{\mathscr T}_{\Omega})$, where $\Gamma:=\{\,||\,.\,||_{\beta,\,l}\,|\,\beta\in\nat^m\,\,\hbox{ and }\,\,l\in\nat\,\}$, is a $G$-locally multiplicatively-convex algebra.

\vskip0.3cm
\small{
\noindent
J. Aragona, R. Fernandez and S. O. Juriaans

\noindent
Instituto de Matem\'atica e Estat\'\i stica

\noindent
Universidade de S\~ao Paulo

\noindent
CP 66281 - CEP 05311-970 - S\~ao Paulo - Brazil

\noindent
aragona@ime.usp.br

\noindent
roselif@ime.usp.br

\noindent
ostanley@ime.usp.br}

\end{document}